\pdfoutput=1 
\documentclass[11pt]{amsart}
\usepackage{epsfig,overpic,comment,mathtools}
\usepackage[usenames,dvipsnames,svgnames,table]{xcolor}
\usepackage[pagebackref,linktocpage=true,colorlinks=true,linkcolor=Blue,citecolor=BrickRed,urlcolor=RoyalBlue]{hyperref}
\usepackage[msc-links,abbrev]{amsrefs} 

\newtheorem{thm}{Theorem}[section]

\newtheorem{problem}[thm]{Problem}

\theoremstyle{definition}

\newtheorem{example}[thm]{Example}
\newtheorem{note}[thm]{Note}

\newcommand{\R}{\mathbf{R}}

\newcommand{\ol}{\overline}
\newcommand{\C}{\mathcal{C}}
\newcommand{\Z}{\mathbf{Z}}

\renewcommand{\S}{\mathbf{S}}
\renewcommand{\l}{\langle}
\renewcommand{\r}{\rangle}

\renewcommand{\tilde}{\widetilde}

\DeclareMathOperator{\Lk}{Lk}
\DeclareMathOperator{\Tw}{Tw}
\DeclareMathOperator{\Wr}{Wr}
\DeclareMathOperator{\Length}{Length}
\DeclareMathOperator{\Cr}{Cr}
\DeclareMathOperator{\SL}{SL}
\DeclareMathOperator{\Rot}{Rot}

\DeclareMathOperator{\sign}{sign}

\title[Topology of closed asymptotic curves]{Topology of closed asymptotic curves\\ on negatively curved surfaces}

\author{Mohammad Ghomi}
\address{School of Mathematics, Georgia Institute of Technology, Atlanta, Georgia 30332}
\email{ghomi@math.gatech.edu}
\urladdr{www.math.gatech.edu/~ghomi}

\author{Matteo Raffaelli}
\address{School of Mathematics, Georgia Institute of Technology, Atlanta, Georgia 30332}
\email{raffaelli@math.gatech.edu}
\urladdr{https://matteoraffaelli.com}

\date{\today \,(Last Typeset)}
\subjclass[2020]{Primary 53A04, 53A05; Secondary 57K10, 35L10}
\date{Last revised on \today}
\keywords{Rigidity of isometric embeddings, Tight surfaces, Characteristic curves,  Hyperbolic PDE, C\u{a}lug\u{a}reanu's formula, Self-linking number, Writhe, Torsion of closed curves.}
\thanks{The first-named author was supported by NSF grant DMS-2202337.}

\begin{document}

\begin{abstract}
Motivated by Nirenberg's problem on isometric rigidity of tight surfaces, we study closed asymptotic curves $\Gamma$ on negatively curved surfaces $M$ in Euclidean $3$-space. In particular, using C\u{a}lug\u{a}reanu's theorem, we obtain a formula for the linking number $\Lk(\Gamma,n)$ of $\Gamma$ with the normal $n$ of $M$. It follows that when $\Lk(\Gamma, n)=0$,  $\Gamma$ cannot have  any locally star-shaped planar projections with vanishing crossing number, which extends observations of  Kovaleva,  Panov  and Arnold. These results hold also for curves with nonvanishing torsion and their binormal vector field. Furthermore
we construct an example where $n$ is injective but $\Lk(\Gamma, n)\neq 0$, and discuss various restrictions on $\Gamma$ when $n$ is injective.
\end{abstract}

\maketitle


\section{Introduction}
An \emph{asymptotic curve} on a negatively curved surface $M$ in Euclidean space $\R^3$ is an integral curve of  the directions where the second fundamental form of $M$ vanishes. These objects form characteristic curves of the hyperbolic PDE for isometric embeddings, and thus play a fundamental role in surface theory. In particular, a well-known problem of Nirenberg \cite{nirenberg1963,han2006,yau2000,rozendorn1992,kantor1980,han-hong-huang2017,han-khuri2011} is concerned with existence of closed asymptotic curves $\Gamma$ when $M$ forms the interior of an annular surface $\ol M$ bounded by convex curves in planes tangent to $\ol M$. Nirenberg conjectured that no such curves exist, which would extend  Alexandrov's rigidity theorem \cite{alexandrov1938} for analytic tight surfaces \cite{cecil-chern1997} to the smooth category. 

To study Nirenberg's problem, we first observe that the Gauss map $n\colon M\to\S^2$ is injective, see Note \ref{note:injective}. 
Furthermore,  the linking number $\Lk(\Gamma,n)=0$ if $M$ is embedded, see Note \ref{note:embedded}. Thus the problem, in the embedded case, may be reduced to the following topological question:

\begin{problem}\label{prob:main}
Let $M\subset\R^3$ be a negatively curved embedded surface with Gauss map $n$ and $\Gamma\subset M$ be a closed asymptotic curve. Can $\Lk(\Gamma, n)= 0$ if $n$ is injective?
\end{problem}

\noindent A negative answer to Problem \ref{prob:main} will immediately prove Nirenberg's conjecture in the embedded case. 
We study Problem \ref{prob:main} using techniques from geometric knot theory, and obtain a number of results which relate the geometry and topology of $\Gamma$ to properties of its planar projections. Some of these results may be of independent interest as they generalize previous observations of Panov \cite{panov1998} and Arnold \cite{arnold1999}, who had noticed some remarkable properties of asymptotic curves without apparent knowledge of Nirenberg's problem. Furthermore, asymptotic curves are closely related to curves with nonvanishing torsion, and our results also apply to these curves (without any reference to a surface). Finally, we will discuss some examples which will further illustrate subtleties of Problem \ref{prob:main}.

To state our main results, let $M$, $n$, and $\Gamma$ be as in Problem \ref{prob:main}. We will assume that $M$ is $\C^3$, which ensures that $\Gamma$ is $\C^2$, see Section \ref{sec:darboux}. 
Let $u\in\S^2$ be a random direction, $\Gamma_u$ be the projection of $\Gamma$ into a plane orthogonal to $u$,   $\Cr(\Gamma_u)$ be  the sum of signed crossings of $\Gamma_u$, $\#\big\{\l u,n\r=0\big\}$ be the number of zeros of $\l u,n\r$, and $\tau_g$ be the geodesic torsion of $\Gamma$, which as we will show in Section \ref{sec:darboux} has a fixed sign. 

\begin{thm}\label{thm:main}
\begin{equation}\label{eq:main}
\Lk(\Gamma,n)=\Cr(\Gamma_u)+\frac{1}{2} \,\#\big\{\l u,n\r=0\big\}\,\textup{sign}(\tau_g).
\end{equation}
\end{thm}

This equation may be regarded as a refinement of C\u{a}lug\u{a}reanu's formula \cite{cualuguareanu1959, cualuguareanu1961,pohl1968,banchoff1976} for the self-linking number of curves with nonvanishing curvature and torsion, see Note \ref{note:self-linking}. We say that $\Gamma_u$ is \emph{locally star-shaped} if its tangent lines do not cover the plane of $\Gamma_u$. Using the above formula we show:

\begin{thm}\label{thm:main2}
If $\Cr(\Gamma_u)=\Lk (\Gamma,n)$, then $\Gamma_u$ cannot be locally star-shaped.
\end{thm}

In particular, if $\Lk(\Gamma, n)=0$, or $M$ is isotopic to a planar annulus, then $\Gamma_u$ cannot be star-shaped, as had been observed by Kovaleva \cite{kovaleva1968}. Panov \cite{panov1998} and Arnold \cite{arnold1999} reproved Kovaleva's result in the case where $M$ forms a graph over the $xy$-plane. Petrunin \cite[p.\ 31]{petrunin2023} also gave a  proof in this case, see Note \ref{note:petrunin}. Our proof of the above theorem leads to our next result. The \emph{absolute rotation index} of a closed curve immersed in the plane is the absolute value of the degree of its Gauss map. An \emph{inflection} is a point of vanishing curvature.

\begin{thm}\label{thm:main3}
If $\Gamma_u$ has no inflections, then it cannot be locally star-shaped; in particular, the absolute  rotation index of $\Gamma_u$ cannot be less than $3$.
\end{thm}

Any space curve with nonvanishing torsion may be realized as an asymptotic curve on a negatively curved surface orthogonal to its binormal vector, see Note \ref{note:darboux}. Thus the above results hold for curves with nonvanishing torsion and their binormal. In particular,
Theorem \ref{thm:main3}  implies that curves of type $(p,1)$ on a torus of revolution cannot have nonvanishing torsion, as had been observed by Costa \cite{costa1990}. This result also generalizes the fact that a curve of nonvanishing torsion cannot project into a strictly convex planar curve, which has been observed  in various contexts \cite{bray-jauregui2015,ivanisvili2023}. See \cite{ghomi2019,ghomi:verticesC,ghomi2017,chen-wang-bin2020} for other recent results on torsion.

When $n$ is injective, $\Gamma$  must  have inflections, see Note \ref{note:inflection}. At these points the Frenet frame of $\Gamma$ is not well-defined, but the Darboux frame of $\Gamma$ with respect to $M$ will be useful, as we describe in Section \ref{sec:darboux}. Applying C\u{a}lug\u{a}reanu's formula to the Darboux frame will yield Theorem \ref{thm:main}, and Theorems \ref{thm:main2} and \ref{thm:main3} follow quickly in Section \ref{sec:proofs}. In Section \ref{sec:examples} we construct an example where $\Lk(\Gamma,n)\neq 0$ but $n$ is injective, which complements constructions by Kovaleva \cite{kovaleva1995} and Arnold \cite{arnold1999} where $\Lk(\Gamma,n)=0$ but $n$ is noninjective. This example also contradicts a claim of Kovaleva \cite{kovaleva1995}, see Note \ref{note:kovaleva}. Other observations concerning Problem \ref{prob:main} will be discussed in Section \ref{sec:injective}.

\begin{note}\label{note:injective}
The Gauss map $n$ in Nirenberg's problem  is injective  on $M$. To see this note that since $M$ has curvature $K<0$ it must lie within the convex hull of $\partial\ol M$. Thus the inward normal  of each component $C_i$ of $\partial \ol M$ points to the interior of the convex body bounded by $C_i$. So the geodesic curvature of $C_i$  in $M$ is the same as the geodesic curvature of $C_i$ in the plane where it lies. Hence, by the convexity  assumption, the total geodesic curvature of $C_i$  is $2\pi$. The Gauss-Bonnet theorem then yields that $\int_{\ol M}K=-4\pi$.
But $n$  is onto, since it maps each $C_i$ to a point, and is locally injective on $M$ since $K\neq 0$.  Hence $n$ must be injective on $M$.
\end{note}

\begin{note}\label{note:embedded}
Let $S$ be the topological sphere obtained by gluing the annular surface $M$ in Nirenberg's problem to the two convex disks bounded by the components of $\partial M$ in their planes. If $M$ is embedded, then there exists $\epsilon>0$ such that the curve $\Gamma_\epsilon:=\Gamma+\epsilon n$ is disjoint from $S$. On the other hand, $\Gamma$ may be contracted to a point in $S$. Hence $\Gamma$ and $\Gamma_\epsilon$ are unlinked, or $\Lk(\Gamma, n)=0$. Here we assume that since $M$ is embedded, $\Gamma$ is simple, and therefore $\Lk(\Gamma, n)$ is well-defined. This assumption is indeed warranted as discussed in Section \ref{sec:reg-simp}.  If $M$ is not embedded, which is possible \cite[p.\ 68]{banchoff-kuhnel1997}, then $\Gamma$ may self-intersect, and so $\Lk(\Gamma, n)$ may not be well-defined.
\end{note}

\section{Preliminaries}\label{sec:darboux}
\subsection{Regularity and simplicity of $\Gamma$}\label{sec:reg-simp}
From now on we assume that $M\subset\R^3$ is a $\C^{k\geq 3}$ embedded orientable surface with Gauss curvature $K<0$, Gauss map $n$,   
second fundamental form $\mathrm{I\!I}(X,Y)\coloneqq \l dn(X),Y\r$, and closed asymptotic curve
$\Gamma\colon\S^1\simeq\R/(2\pi\Z)\to M$. We will write $\Gamma$ to refer both to the mapping and its image $\Gamma(\S^1)$, and will write $n(t)$ to mean $n(\Gamma(t))$.
We assume that the speed $|\Gamma'|=1$, after a rescaling, and set $T\coloneqq \Gamma'$. Then
$$
\mathrm{I\!I}(T,T)=0.
$$
This equation may be written as $\ell_{11}(du_1/du_2)^2+2\ell_{12}(du_1/du_2)+\ell_{22}=0$ in local coordinates where $(u_1,u_2)$ represent $\Gamma$ and $\ell_{ij}$ are the coefficients of $\mathrm{I\!I}$. Thus  $\Gamma$ is $\C^{k-1}$ since 
 $\ell_{ij}$ are $\C^{k-2}$. It is easy to see that one can distinguish precisely two fields of asymptotic lines on $M$ \cite[p.\ 184]{nirenberg1963}. Thus, by the uniqueness of solutions to ODE, $\Gamma$ cannot self-intersect in $M$. Since $M$ is embedded in $\R^3$, it follows that $\Gamma$ is simple.
 
 \subsection{The Darboux frame}
 Set $n^\perp\coloneqq n\times T$. Then
$(T,n^\perp,n)$ forms the \emph{Darboux frame} of $\Gamma$ with respect to $M$. Since $M$ is $\C^3$ and $\Gamma$ is $\C^2$, $(T,n^\perp,n)$ is $\C^1$, and we have
\begin{equation}\label{eq:darboux}
T'=\kappa_g n^\perp,\quad\quad (n^\perp)'=-\kappa_g T+\tau_g n,\quad\quad n'=-\tau_g n^\perp,
\end{equation}
where $\kappa_g\coloneqq  \l T', n^\perp\r$ and $\tau_g\coloneqq \l(n^\perp)',n\r$ are the \emph{geodesic curvature} and \emph{geodesic torsion} of $\Gamma$ respectively. Let $\kappa\coloneqq |T'|$ be the curvature of $\Gamma$. Then $|\kappa_g|=\kappa$. So there is no distinction between inflections of $\Gamma$ with respect to $\kappa$ or $\kappa_g$. If $\kappa\neq 0$ at any point, then the principal normal $N\coloneqq T'/ \kappa$ and the binormal $B\coloneqq T\times N$ of $\Gamma$ generate the \emph{Frenet frame} $(T,N,B)$.
We have $(N,B)=\pm(n^\perp, n)$ depending on the sign of $\kappa_g$. So when $\kappa\neq 0$ the Darboux and Frenet frames coincide up to reflections.
In particular, when $k\geq4$, or $\Gamma$ is $\C^3$, the torsion $\tau\coloneqq \l N',B\r$ of $\Gamma$ coincides with $\tau_g$. 
Note that $|\tau_g|=|n'|=|dn(T)|$, and $dn$ is nondegenerate since $K\neq 0$. Thus
\begin{equation}\label{eq:tau}
\tau_g\neq 0.
\end{equation}
This observation also follows from the Beltrami-Enneper theorem \cite[p.\ 200]{spivak:v3} \cite[p.\ 609]{hartman-wintner1952}. 

\subsection{Crossings}
For any simple closed immersion $\Gamma\colon\S^1\to\R^3$ and direction $u\in\S^{2}$, let $\Gamma_u$ denote the projection of $\Gamma$ into a plane $\Pi$ orthogonal to $u$. For almost every $u$, or \emph{random direction}, $\Gamma_u$ is in \emph{general position}. More explicitly, there are only finitely many points $p_i\in\Pi$, called \emph{crossings}, such that $\Gamma_u^{-1}(p_i)$ consists of more than one point. Furthermore all $p_i$ will be transversal double points, i.e., $\Gamma_u^{-1}(p_i)=\{t_i^+,t_i^-\}$ with $\Gamma_u'(t_i^+)\times \Gamma_u'(t_i^+)\neq 0$. We assume that $\l \Gamma(t_i^+),u\r>\l \Gamma(t_i^-),u\r$. Then $\sign(p_i)$ is defined as the sign of $\l\Gamma_u(t_i^+)\times\Gamma_u(t_i^-),u\r$, i.e., $1$ or $-1$ depending on whether  $\l\Gamma_u(t_i^+)\times\Gamma_u(t_i^-),u\r>0$ or $<0$ respectively. Then the \emph{crossing number} of $\Gamma_u$ is given by $\Cr(\Gamma_u)\coloneqq \sum\sign(p_i)$.

For any pair of disjoint immersions $\Gamma^1$, $\Gamma^2\colon\S^1\to\R^3$, the \emph{crossing number} $\Cr(\Gamma^1_u,\Gamma^2_u)$ is defined similarly. Again, assuming $u$ is a random direction or $\Gamma^1$ and $\Gamma^2$ are in general position, then $\Gamma_u^1\cap\Gamma_u^2$ consists of only a finite number of points $p_i$, and for each $p_i$ there will be exactly one pair of points $t_i^1$, $t_i^2\in\S^1$ such that $\Gamma_u^1(t_i^1)=\Gamma_u^2(t_i^2)=p_i$. Then $\sign(p_i)$ is defined as the sign of $\l(\Gamma_u^1)'(t_i^1)\times(\Gamma_u^2)'(t_i^2),u\r$ or $\l(\Gamma_u^2)'(t_i^2)\times\l(\Gamma_u^1)'(t_i^1),u\r$ if $\l \Gamma^1(t_i^1),u\r>\l \Gamma^2(t_i^2),u\r$ or $\l \Gamma^1(t_i^1),u\r<\l \Gamma^2(t_i^2),u\r$ respectively. Finally we set $\Cr(\Gamma^1_u,\Gamma^2_u)\coloneqq \sum\sign(p_i)$.

Let $v$ be a unit normal vector field along $\Gamma$. The pair $(\Gamma, v)$ is called a \emph{ribbon} based on $\Gamma$. Let $\epsilon>0$ be so small that the perturbation $\Gamma+\epsilon v$ is disjoint from $\Gamma$. The crossings $\Gamma_u\cap (\Gamma +\epsilon v)_u$ fall into two categories \cite{dennis-hannay2005}: a crossing formed at $t\in\S^1$ is \emph{local} if $\Gamma_u(t)=(\Gamma(t)+\epsilon v(t))_u$, or $v(t)=\pm u$; otherwise, it is  \emph{nonlocal}. Nonlocal crossings converge in pairs to self-crossings of $\Gamma_u$ as $\epsilon\to 0$. Thus 
\begin{equation}\label{eq:crossings}
\Cr\big(\Gamma_u, (\Gamma+\epsilon v)_u\big)=2\Cr(\Gamma_u)+\Cr_{local}\big(\Gamma_u, (\Gamma +\epsilon v)_u\big),
\end{equation}
where $\Cr_{local}$  is the signed sum of local crossings.

\subsection{C\u{a}lug\u{a}reanu's formula}
For a ribbon $(\Gamma, v)$, 
C\u{a}lug\u{a}reanu's formula \cite{moffatt-ricca1992,dennis-hannay2005,fuller1971,aew2004} states that
$$
\Lk(\Gamma, v)= \Wr(\Gamma)+ \Tw(\Gamma,v),
$$
where $\Lk$, $\Wr$ and $\Tw$ stand for the linking number, writhe and twist respectively. The linking number is defined as
$
\Cr\big(\Gamma_u, (\Gamma +\epsilon v)_u\big)/2.
$
So by \eqref{eq:crossings}
\begin{equation}\label{eq:link1}
\Lk(\Gamma, v)=\Cr(\Gamma_u)+\frac{1}{2}\Cr_{local}\big(\Gamma_u, (\Gamma +\epsilon v)_u\big).
\end{equation}
 The \emph{twist} of $(\Gamma,v)$ is given by 
$
\int_\Gamma \l (v^\perp)',v\r/(2\pi)
$
where $v^\perp\coloneqq v\times T$. Let $\theta_v\colon\S^1\to\R$  be the continuous function such that
$
v(t)=\cos(\theta_v(t))n^\perp(t)+\sin(\theta_v(t))n(t),
$  
and
$
\Rot(v,n)\coloneqq (\theta_v(L)-\theta_v(0))/(2\pi)
$
denote the \emph{total rotation} of $v$ with respect to $n$. 
A computation using \eqref{eq:darboux} yields that 
$$
\Tw(\Gamma,v)
=\frac{1}{2\pi}\int_\Gamma \tau_g+ \Rot(v,n).
$$
Finally, the \emph{writhe} of $\Gamma$ is 
the average of the crossing numbers of $\Gamma_u$,
\begin{equation}\label{eq:writhe}
\Wr(\Gamma)\coloneqq \frac{1}{4\pi}\int_{u\in\S^2}\Cr(\Gamma_u).
\end{equation}
Note that 
$
\Wr(\Gamma)=\Lk(\Gamma, n)-\Tw(\Gamma,n)=\Lk(\Gamma, n)- \int_\Gamma\tau_g/(2\pi).
$
Thus we conclude that
\begin{equation}\label{eq:link2}
 \Lk(\Gamma, v)= \Lk(\Gamma,n)+ \Rot(v,n).
\end{equation}

\begin{note}\label{note:darboux}
If a curve $\Gamma$ in $\R^3$ admits an orthonormal frame satisfying \eqref{eq:darboux} and \eqref{eq:tau}, e.g., $\Gamma$ has nonvanishing curvature and torsion, then  a routine computation shows that $\Gamma(t)+s n^\perp(t)$ generates a negatively curved surface, for $-\epsilon <s<\epsilon$, which contains $\Gamma$ as an asymptotic curve. More specifically, $K(t,0)=-\tau_g^2(t)$. Thus any curve satisfying \eqref{eq:darboux} and \eqref{eq:tau} may be called an asymptotic curve (without reference to a specific surface), and any curve with nonvanishing curvature and torsion is an asymptotic curve.
\end{note}

\section{Proofs of the Main Results}\label{sec:proofs}
\subsection{Proof of Theorem \ref{thm:main}}
Let $u\in\S^2$ be a direction which is not parallel to any tangent line of $\Gamma$, and set $u^\perp\coloneqq  u\times T/|u\times T|$. By \eqref{eq:link2}
 $$
\Lk(\Gamma,n)= \Lk(\Gamma, u^\perp)-\Rot(u^\perp,n).
 $$
Note that $(\Gamma+\epsilon u^\perp)_u$ is locally disjoint from $\Gamma_u$. Thus  \eqref{eq:link1} yields that
 $$
 \Lk(\Gamma, u^\perp)=\Cr(\Gamma_u).
 $$
It remains then to compute $\Rot(u^\perp,n)$. To this end consider the mapping $\nu\colon\S^1\to\S^1$ given by $\nu(t)\coloneqq e^{i\theta(t)}$, where $\theta\coloneqq \theta_{u^\perp}$ is as defined above, i.e., $u^\perp(t)=\cos(\theta(t))n^\perp+\sin(\theta(t))n$. Then
  $
\Rot(u^\perp,n)=\deg(\nu).
$
To compute  $\deg(\nu)$  assume that $\S^1$ is oriented counterclockwise. So $it_0\in T_{t_0}\S^1$ has positive orientation. We have $d\nu_{t_0}(it_0)=i\nu(t_0)\theta'(t_0)\in T_{\nu(t_0)}\S^1$. To find $\theta'$ note that
 $$
 \cot(\theta)
 =
 \frac{\l u^\perp, n^\perp\r}{\l u^\perp,  n\r}
 =
 \frac{\l u\times T, n^\perp\r}{\l u\times T, n\r}
 =
 \frac{\l u,T\times n^\perp\r}{\l u, T\times n\r}
 =
- \frac{\l u,n \r}{\l u,n^\perp\r}
.
 $$
Differentiating the far sides of this expression using \eqref{eq:darboux}, we obtain
\begin{gather*}
\frac{\theta'}{\sin^2(\theta)}
=
-\tau_g-\frac{\l u, n\r\l u, -\kappa_gT+\tau_gn^\perp\r}{\l u,n^\perp\r^2}.
\end{gather*}
Suppose $\nu(t_0)=(0,\pm1)$, or $\theta(t_0)= m\pi+\pi/2$, for $m\in\textbf{Z}$. Then $u^\perp(t_0)=n(t_0)$, which yields $\l u,n(t_0)\r=0$. It follows that $\theta'(t_0)=-\tau_g(t_0)$. So $d\nu_{t_0}(it_0)=-i\nu(t_0)\tau_g(t_0)$, which does not vanish by \eqref{eq:tau}.
Thus  $(0,\pm 1)$ are regular values of $\nu$, and $d\nu_{t_0}$ preserves orientation if and only if $\tau_g(t_0)<0$. Hence, since $\tau_g$ has constant sign, 
$
-\textup{sign}(\tau_g)\deg(\nu)
=\#\{\nu^{-1}(0,\pm1)\}/2= \#\big\{\l u,n\r=0\big\}/2.
$
So
$$
\Rot(u^\perp,n)=-\frac{1}{2}\, \#\big\{\l u,n\r=0\big\}\,\textup{sign}(\tau_g)
$$
which completes the proof. 

\subsection{Proof of Theorem \ref{thm:main2}}
If $\Cr(\Gamma_u)=\Lk(\Gamma,n)$ then $\langle u,n\rangle\neq 0$ by \eqref{eq:main}. So $M$ is locally a graph over a plane orthogonal to $u$, which we may identify with the $xy$-plane. After an affine transformation, given by $(x,y,z)\mapsto (x,y,\lambda z)$, we may assume that $M$ is arbitrarily close to the $xy$-plane, since affine transformations preserve asymptotic curves. Suppose, towards a contradiction, that $\Gamma_u$ is locally star-shaped with respect to the origin, i.e., tangent lines of $\Gamma_u$ do not pass through $o$. Then $\l \Gamma_u, u^\perp\r\neq 0$. Let $h\coloneqq \langle \Gamma,n\rangle$. Then $h'=-\tau_g\langle \Gamma, n^\perp\rangle$. But $\l\Gamma,n^\perp\r\to\l\Gamma_u,u^\perp\r$ as $\lambda\to 0$. So $h'\neq 0$, for $\lambda$ small, which is impossible since $\Gamma$ is closed.

\subsection{Proof of Theorem \ref{thm:main3}} We may again assume that $u=(0,0,1)$ and $\Gamma_u$ lies in the $xy$-plane.  When the curvature of $\Gamma_u$ does not vanish, the curvature $\kappa$ of $\Gamma$ does not vanish either, and the principal normal $N\neq \pm u$.  But if $\kappa$ does not vanish, then $n^\perp=\pm N$, as discussed in the last section. So $n^\perp\neq\pm u$. Then, after  the rescaling $(x,y,z)\mapsto (x,y,\lambda z)$, we may assume that $|\l T, u\r|$ and $|\l n^\perp, u\r|$ are arbitrarily small, which  implies that $n$ is almost parallel to $u$ or $-u$. In particular $\l n, u\r\neq 0$. Thus $M$ is locally a graph over the $xy$-plane. Now if $\tau_g\neq 0$, then we again arrive at a contradiction, as shown in the proof of Theorem \ref{thm:main2}.

\begin{note}\label{note:petrunin}
An alternative argument for finishing the proofs of Theorems \ref{thm:main2} and \ref{thm:main3}, once we know that $M$ is locally a graph over the $xy$-plane, can be given following  Petrunin \cite[p.~31]{petrunin2023}. Let $u=(0,0,1)$ and set
$
h\coloneqq \l\Gamma,n\r/\l u,n \r,
$
i.e., the height over the origin of the tangent planes of $M$ along $\Gamma$. We compute that
$$
-h' \l u,n\r^2
=
\tau_g\big(\l\Gamma,n^\perp\r \l u,n \r-\l\Gamma,n\r\l u,n^\perp \r\big)
=
\tau_g\l\Gamma,T\times u\r
=
\tau_g|T\times u|\l\Gamma,u^\perp\r.
$$
If $\Gamma_u$ is locally star-shaped with respect to the origin, then
$
\l\Gamma,u^\perp\r 
=
\l\Gamma_u,u^\perp\r
\neq 0.
$
So $h'\neq 0$, which is again a  contradiction.
\end{note}

\begin{note}\label{note:self-linking}
As we mentioned in Section \ref{sec:darboux}, if $\Gamma$ has no inflections, then $n=\pm B$ and $\tau_g=\tau$. So $\Lk(\Gamma,n)=\Lk(\Gamma, B)=\Lk(\Gamma,N)$ which is known as the self-linking number of $\Gamma$ and is denoted by $\SL(\Gamma)$ \cite{pohl1968,rogen2008}. Thus \eqref{eq:main} yields that
\begin{equation}\label{eq:SL}
\SL(\Gamma)=\Cr(\Gamma_u)+\frac{1}{2} \#\big\{\l u,B\r=0\big\}\sign(\tau).
\end{equation}
By Crofton's formula, 
$$
\int_{u\in\S^2}  \#\big\{\l u,B\r=0\big\}\sign(\tau)= 2\Length(B)\sign(\tau)=2\int_\Gamma\tau.
$$
Thus \eqref{eq:SL} together with \eqref{eq:writhe} yields
$
\SL(\Gamma)=\Wr(\Gamma)+ \int_\Gamma\tau/(2\pi),
$
which is a special case of C\u{a}lug\u{a}reanu's formula \cite{pohl1968}. Hence \eqref{eq:main} generalizes C\u{a}lug\u{a}reanu's formula for self-linking number of curves with positive curvature and torsion. For this class of curves \eqref{eq:SL} may be rewritten as
\begin{equation}\label{eq:banchoff}
\SL(\Gamma)=\Cr(\Gamma_u)+\frac{1}{2} \textup{Inflection}(\Gamma_u),
\end{equation}
where $\textup{Inflection}(\Gamma_u)$ denotes the number of inflections of $\Gamma_u$.
Indeed $\Gamma_u(t)$ is an inflection if and only if $N(t)$ projects into the tangent line of $\Gamma_u$ at $\Gamma_u(t)$, or $0=\l N(t), u^\perp(t)\r=|u\times T(t)|\l u, B(t)\r$. Formula \eqref{eq:banchoff} had been observed by Banchoff \cite[p.\ 1182]{banchoff1976}, and also follows from Milnor \cite{milnor1953} who had shown that 
$\int_\Gamma|\tau|$
is the average value of $\textup{Inflection}(\Gamma_u)$, when $\Gamma$ has no inflections.
\end{note}

\section{Examples}\label{sec:examples}
Here we describe a pair of examples which illustrate some of the subtleties of Problem \ref{prob:main}.
\begin{example}
For the sake of comparison we start with Kovaleva's example \cite{kovaleva1968},  which is given by the coordinate functions
\begin{eqnarray*}
\Gamma_1(t)&\coloneqq &(3+\sin(t))\cos(\sqrt{63/8}\cos(t)),\\
\Gamma_2(t)&\coloneqq &(3+\sin(t))\sin(\sqrt{63/8}\cos(t)), \\
\Gamma_3(t)&\coloneqq &\sin(2t)+46\cos(t)-27\cos^3(t)+27/8\cos^5(t),
\end{eqnarray*}
see Figure \ref{fig:kovaleva}(a).
\begin{figure}[h]
\centering
 \begin{overpic}[height=1.3in]{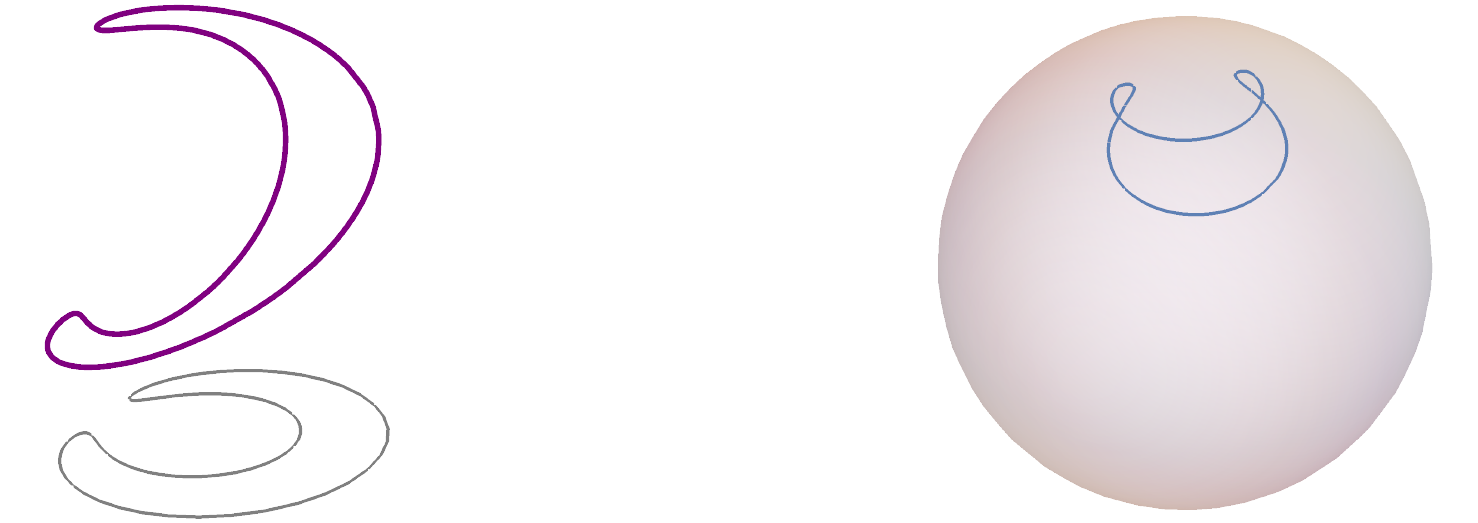}
 \put(12,-4){(a)}
 \put(78,-4){(b)}
 \end{overpic}
  \caption{}\label{fig:kovaleva}
\end{figure}
We have provided a Mathematica notebook \cite{ghomi-raffaelli-mathematica2024} where one can check that this curve admits a Darboux framing $(T,n^\perp,n)$ with nonvanishing torsion $\tau_g$, as described in Note \ref{note:darboux}, and therefore is an asymptotic curve.  Figure 
\ref{fig:kovaleva}(b) shows that $n(\Gamma)$ is not injective. Furthermore, by \eqref{eq:main} we have $\Lk(\Gamma, n)=0$, since as Figure \ref{fig:kovaleva} shows,
$\Cr(\Gamma_u)=0=\#\{\l n,u\r=0\}$ for $u=(0,0,1)$.

Aicardi \cite{aicardi2000} showed that curves of nonvanishing curvature and torsion can be constructed with any self-linking number. Thus there exists closed asymptotic curves with any self liking number, as we discussed in Note \ref{note:darboux}. The special feature of Kovaleva's example and similar constructions by Arnold, however, is that the surface forms a graph, which makes these examples much more subtle.

\end{example}

\begin{example}\label{ex:2}
Next we construct a closed asymptotic curve with opposite properties, i.e., $n$  injective but $\Lk(\Gamma, n)\neq 0$. Again the reader can verify our computations via the Mathematica notebook \cite{ghomi-raffaelli-mathematica2024} that we have provided. This example is obtained by starting with an embedding $n\colon \S^1\to\S^2$ given by 
\begin{eqnarray*}
n_1(t)&\coloneqq&(3+\sin(t))\cos(5/2\cos(t)),\\
n_2(t)&\coloneqq &(3+\sin(t))\sin(5/2\cos(t)), \\
n_3(t)&\coloneqq &(1-n_1^2(t)-n_2^2(t))^{1/2},
\end{eqnarray*}
see  Figure \ref{fig:rocket}(a).
\begin{figure}[h]
\centering
 \begin{overpic}[height=1.5in]{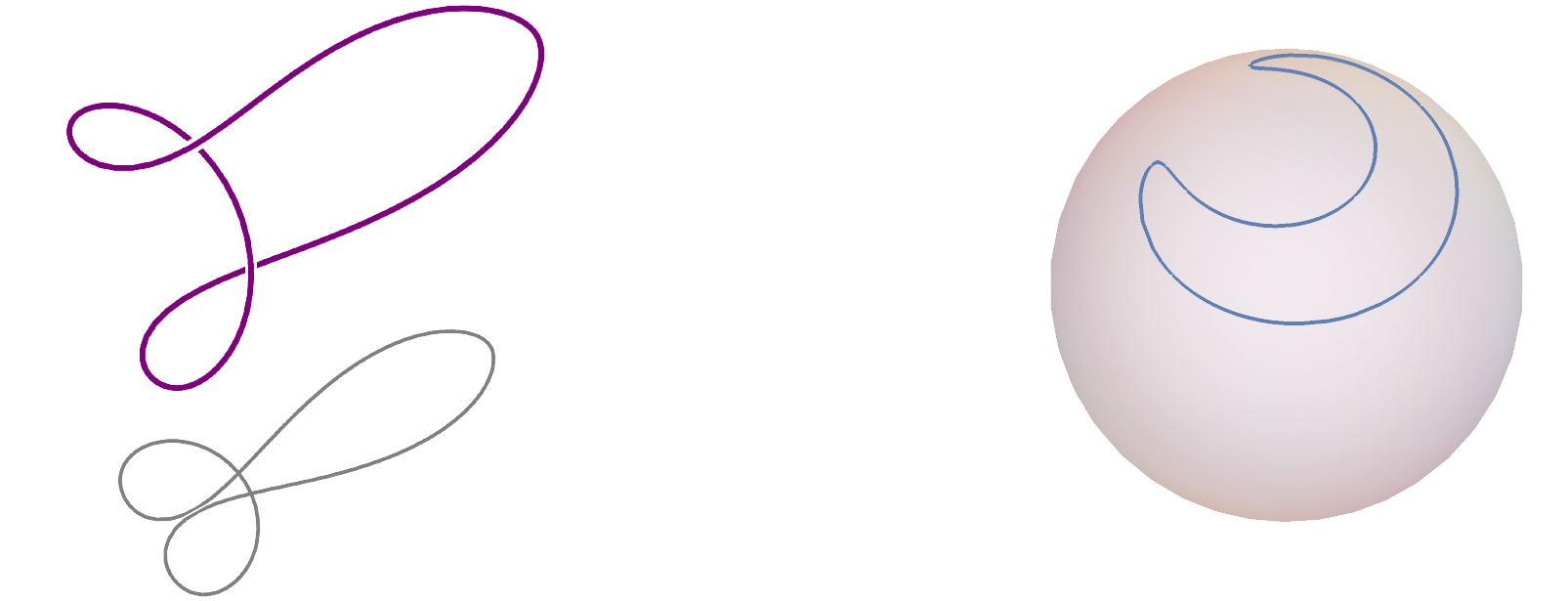}
\put(17,0){(a)}
 \put(80,0){(b)}
 \end{overpic}
  \caption{}\label{fig:rocket}
\end{figure}
Set $\tau_g\coloneqq |n'|$, and $T\coloneqq -n\times n'/\tau_g$. Note that $T$ traces the center of oriented great circles  
tangent to $n$. We have constructed $n$ so that these circles cover $\S^2$, i.e., $n$ is not star-shaped. It follows that the origin is contained in the interior of the convex hull of $T$, see Figure \ref{fig:T}.
 Thus there exists a $\C^\infty$ positive  function $\rho\colon \S^1\to\R$ with $\int\rho T=0$, which can be constructed using methods of convex integration \cite[Lem.\ 2.3]{ghomi2007} \cite[p.\ 168]{gromov1986}. Setting $\Gamma(t)\coloneqq \int_0^{t}\rho(s)T(s)ds$, for an appropriate choice of $\rho$, yields the desired curve depicted in Figure \ref{fig:rocket}(a). Here $\Cr(\Gamma_u)=2$ whereas, $\#\{\langle n, u\rangle=0\}=0$. Thus by \eqref{eq:main}, $\Lk(\Gamma,n)=2$. 
 
 It remains to find $\rho$. Let $t_1=7\pi/6$, $t_2=11\pi/6$, $t_3=\pi/6$, $t_4=5\pi/6$, and $t_5=\pi/2$. Then the origin $o$ of $\R^3$ is contained in the interior of the convex hull of $p_i\coloneqq T(t_i)$, see Figure \ref{fig:T}.    
 \begin{figure}[h]
\centering
 \begin{overpic}[height=1.2in]{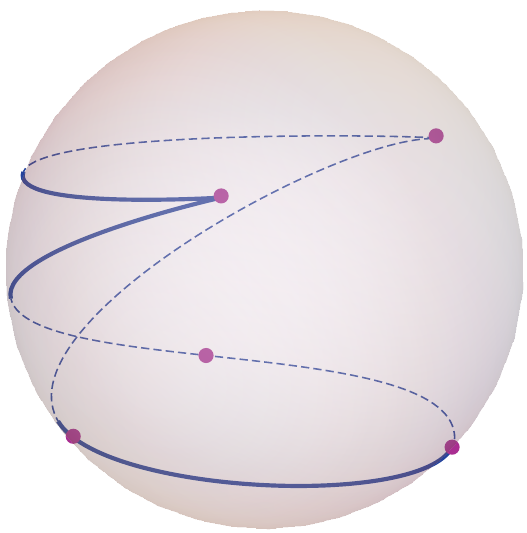}
 \put(79,79){\small $p_1$}
\put(38,67){\small $p_2$}
\put(35,37){\small $p_3$}
\put(6,22){\small $p_4$}
\put(87,10){\small $p_5$}
 \end{overpic}
  \caption{}\label{fig:T}
\end{figure}
One may also note that the pairs $(p_1,p_2)$ and $(p_3,p_4)$ are symmetric with respect to the $yz$-pane, while $p_5$ lies on the $yz$-plane. Furthermore, $p_1$ and $p_2$ lie on the cusps of $T$ which correspond to the inflections of $n$.
 Set $p_i'\coloneqq \int\phi_i T$ for positive functions $\phi_i\colon \S^1\to\R$ with $\int\phi_i=1$. If $\phi_i$ are concentrated near $t_i$ then $|p_i'- p_i|$ are small. So $o$ lies in the interior of the convex hull of $p_i'$. For instance we may set $\phi_i(t)\coloneqq (1.1+\cos(t-t_i))^{10}/\int(1.1+\cos(t-t_i))^{10}$. Then  there exist constants $c_i>0$ such that $\sum c_i p_i'=o$. In particular we may set $c_1=c_2=1$, $c_3=c_4\approx 0.0815$, and $c_5\approx 0.7465$. Then we set $\rho\coloneqq \sum c_i\phi_i$.
\end{example}

\section{Further Observations}\label{sec:injective}
Here we include some observations on the structure of $\Gamma$ and its projections $\Gamma_u$ in the case where $n$ is injective.
\begin{note}\label{note:inflection}
When $n$ is injective, $\Gamma$ must have inflections. To see this note that the geodesic curvature of $n$ in $\S^2$ is given by 
\begin{equation}\label{eq:geodesic-curvature}
\tilde\kappa_g
\coloneqq  
\frac{\l n'',n\times n'\r}{|n'|^3}
=
\frac{\l -\tau'_gn^\perp-\tau_g (n^\perp)',\tau_gT\r}{\tau_g^3}
=
\frac{\l -(n^\perp)',T\r}{\tau_g}
=
\frac{\kappa_g}{\tau_g},
\end{equation}
and recall that $|\kappa_g|=\kappa$. Thus if $\kappa\neq 0$, then $\kappa_g$ has a fixed sign. Consequently $n$  bounds a geodesically convex domain $\Omega$ in $\S^2$. In particular $n$ lies in a hemisphere centered at a point of $\Omega$. It follows that $T$ lies in the opposite hemisphere, since $T$ is traced by the centers of oriented great circles tangent to $n$. But since $\Gamma$ is closed, the origin of $\R^3$ must lie in the relative interior of the convex hull of $T$. Hence $T$ must be a great circle, which in turn implies that $\Gamma$ is a planar curve, or $n$ is constant, which is a contradiction.
\end{note}

\begin{note}
An alternative argument for showing that injectivity of $n$ forces inflections along $\Gamma$ has been given by Kovaleva \cite{kovaleva1995}. Since that work has been published only in Russian, we include the argument here. 
Let $\tilde\Gamma\coloneqq n(\Gamma)$, and $A_1$, $A_2$ be the areas of the components of $\S^2\setminus\tilde\Gamma$. Then the Gauss-Bonnet theorem yields that
\begin{equation}\label{eq:gauss-bonnet}
\left|\int_{\tilde\Gamma} \tilde\kappa_g\right| =\frac{1}{2}|A_1-A_2|< 2\pi.
\end{equation}
Furthermore, since $|n'|=\tau_g$ it follows from \eqref{eq:geodesic-curvature} that
\begin{equation}\label{eq:kappa-kappa}
\int_{\tilde\Gamma}\tilde\kappa_g
=
\int_0^{2\pi}\tilde\kappa_g(t)|n'(t)|dt
=
\sign(\tau_g)\int_0^{2\pi}\frac{\kappa_g(t)}{\tau_g(t)}\tau_g(t)dt
=
\sign(\tau_g)\int _\Gamma\kappa_g.
\end{equation}
So $|\int\kappa_g|< 2\pi$. Now suppose that $\kappa_g$ does not change sign. Then
$\int \kappa=\int|\kappa_g|=|\int\kappa_g|< 2\pi$. But since $\Gamma$ is a closed curve in $\R^3$,  $\int\kappa\geq 2\pi$ by Fenchel's theorem. Hence  we arrive at a contradiction.
\end{note}

\begin{note}
When $n$ is injective, 
\begin{equation}\label{eq:kappa-tau}
\left(\int_\Gamma \kappa_g\right)^2+\left(\int_\Gamma \tau_g\right)^2>4\pi^2.
\end{equation}
To see this let $L$ denote the length of $\tilde\Gamma=n(\Gamma)$, and $A$ be the area of the region bounded by $\tilde\Gamma$ into which $n\times n'=-\tau_g n\times n^\perp=\tau_g T$ points. Assuming $\tau_g>0$, we have $L=\int_\Gamma|n'|=\int_\Gamma\tau_g$. Furthermore, by Gauss-Bonnet theorem and \eqref{eq:kappa-kappa} $A=2\pi-\int_{\tilde\Gamma} \tilde\kappa_g=2\pi-\int_\Gamma\kappa_g$. By the isoperimetric inequality on $\S^2$,
$L^2\geq 4\pi A -A^2$ with equality only if $\tilde\Gamma$ is a circle.  Since $\tilde\Gamma$ has inflections, as discussed in Note \ref{note:inflection}, $\tilde\Gamma$ cannot be a circle.
Thus 
$L^2> 4\pi A -A^2$, which yields \eqref{eq:kappa-tau}.
\end{note}

\begin{note}
If $n$ is injective and $\langle u,n\rangle\neq 0$, or $M$ is locally a graph over a plane orthogonal to $u$, then the absolute rotation index of $\Gamma_u$ is $1$ (which was precisely the case in Example \ref{ex:2}). To see this suppose again that $u=(0,0,1)$ and consider the rescaling $(x,y,z)\mapsto(x,y,\lambda z)$ as in the proofs of Theorems \ref{thm:main2} and \ref{thm:main3}. Then, as $\lambda\to 0$, $\int \kappa_g$ converges to the total geodesic curvature of $\Gamma_u$. But $\int\kappa_g=\pm\int\tilde\kappa_g$ by \eqref{eq:kappa-kappa}. Furthermore, $\int\tilde\kappa_g$ converges to $\pm 2\pi$ by \eqref{eq:gauss-bonnet}, since $\tilde\Gamma$ converges to $u$ or $-u$. So the total geodesic curvature of $\Gamma_u$ is $\pm 2\pi$, which means the rotation index is $\pm 1$.
\end{note}

\begin{note}\label{note:kovaleva}
In \cite{kovaleva1995} Kovaleva studied Problem \ref{prob:main} in the case where $\Gamma$ has only finitely many inflections, and claimed  that when $n$ is injective, $\Lk(\Gamma,n)$ is equal to half the number of  inflections. But as we pointed out in Note \ref{note:inflection}, inflections of $\Gamma$ correspond to inflections of $n$ in $\S^2$. As Figure \ref{fig:rocket}(b) shows, $n$ has only two inflections in Example \ref{ex:2}, whereas $\Lk(\Gamma, n)=2$ as we discussed above. Thus Example \ref{ex:2} contradicts Kovaleva's claim.
\end{note}

\bibliographystyle{abbrv}
\bibliography{references}

\end{document}